\documentclass[12pt]{article}
\usepackage{amsfonts,amssymb}
\usepackage[mathscr]{eucal}
\usepackage{amsbsy}

\makeindex \voffset = -18mm
\begin{document}
\baselineskip = 5mm
%\tableofcontents
%%%%%%%%%%%%% theorems, lemmas, etc %%%%%%%%%%%%%%%%%%%%%%%%%%%%%%%%%%%%%%%
\newtheorem{theorem}{Theorem}
\newtheorem{lemma}[theorem]{Lemma}
\newtheorem{corollary}[theorem]{Corollary}
\newtheorem{example}[theorem]{Example}
\newtheorem{proposition}[theorem]{Proposition}
\newtheorem{remark}[theorem]{Remark}
\newtheorem{definition}[theorem]{Definition}
\newtheorem{conjecture}[theorem]{Conjecture}
%%%%%%%%%%%%%%%%%%%%%%%%%%%%%%%%%%%%%%%%%%%%%%%%%%%%%%%%%%%%%%%%%%%%%%%%%%%
\newcommand \pal {\! \mid \! }
%%%%%%%%%%%%%%%%%%%%%%%%%%%%%%%%%%%%%%%%%%%%%%%%%%%%%%%%%%%%%%%%%%%%%%%%%%%
\newenvironment{pf}{\par\noindent{\em Proof}.}{\hfill\framebox(6,6)
\par\medskip}
%%%%%%%%%%%%%%%%%%%%%%%%%%%%%%%%%%%%%%%%%%%%%%%%%%%%%%%%%%%%%%%%%%%%%%%%%%%

\title{\bf Finite dimensional motives \\ and the Conjectures of Beilinson
and Murre}

\author{Vladimir Guletski\v \i \hspace{1mm}and Claudio Pedrini
\thanks{Supported by TMR ERB FMRX CT-97-0107 and
INTAS-99-00817. The second named author is a member of GNSAGA of
CNR.}}

\maketitle

\section{Introduction}
\label{intro}

Let $k$ be a field of characteristic $0$ and let $\mathcal V_k$ be
the category of smooth projective varieties over $k$. By $\sim $
we denote an {\it adequate} equivalence relation for algebraic
cycles on varieties \cite{J2}. For every $X\in \mathcal V_k$ let
$A^i_{\sim }(X)=(Z^i(X)/ \sim )\otimes {\mathbb Q}$ be the Chow
group of codimension $i$ cycles on $X$ modulo the chosen relation
$\sim $ with coefficients in $\mathbb Q$.

Let $X,Y\in \mathcal V_k$, let $X=\cup X_i$ be the connected
components of $X$ and let $d_i=\dim (X_i)$. Then $Corr^r_{\sim
}(X,Y)=\oplus _iA^{d_i+r}_{\sim }(X_i\times Y)$ is called a space
of correspondences of degree $r$ from $X$ into $Y$. For any $f\in
Corr^r(X,Y)$ and $g\in Corr^s(Y,Z)$ their composition $g\circ f\in
Corr^{r+s}(X,Z)$ is defined by the formula $g\circ
f={p_{XZ}}_*(p_{XY}^*(f)\cdot p_{YZ}^*(g))$ where $p_{XZ}$,
$p_{XY}$ and $p_{YZ}$ are the appropriate projections. In particular,
we have a linear action of correspondences $Corr^s(Y,Z)\times
A^t(Y)\rightarrow A^{s+t}(Z)$ defined by the rule $(\alpha
,x)\mapsto {p_Y}_*(\alpha \cdot p_X^*(x))$, where $p_X$ and $p_Y$
are the projections.

The category of {\it pure motives} $\mathcal M_{\sim }$ over $k$
with coefficients in $\mathbb Q$ with respect to the given
equivalence relation $\sim $ can be defined as follows \cite{Sch}.
Its objects are triples $M=(X,p,m)$, where $X\in \mathcal V_k$,
$p\in Corr^0_{\sim }(X,X)$ is a projector (i.e. $p\circ p=p$) and
$m\in \mathbb Z$. Morphisms from $M=(X,p,m)$ into $N=(Y,q,n)$ in
$\mathcal M_{\sim }$ are given by correspondences $f\in
Corr^{n-m}_{\sim }(X,Y)$, such that $f\circ p=q\circ f=f$, and
compositions of morphisms are induced by compositions of
correspondences.

The category $\mathcal M_{\sim }$ is pseudoabelian and $\mathbb
Q$-linear. Moreover, it is a tensor category with tensor structure
defined by the formula $(X,p,m)\otimes (Y,q,n)=(X\times Y,p\otimes
q,m+n)$. The triple $1=(Spec(k),id,0)$ plays a role of the unite
object in $\mathcal M_{\sim }$ and the Lefschetz motive $\mathbb
L$ is the triple $(Spec(k),id,-1)$. For any motive $M=(X,p,m)$ one
defines the Tate twist $M(r)$ to be the motive $M\otimes \mathbb
L^{-r}=(X,p,m+r)$, where $\mathbb L^r=\mathbb L^{{\otimes }^r}$
for a positive integer $r$, $\mathbb L^0=1$ and $\mathbb
L^r=\mathbb L^{{\otimes }^{-r}}$ for a negative $r$, see
\cite{Sch}, 1.9. At last, $\mathcal M_{\sim }$ is rigid \cite{J2}
in the sense that there exists internal $Hom$'s and dual obiects
$M^*$ for all $M\in \mathcal M_{\sim }$ satisfying well known
axioms \cite{DeMi}.

For any algebraic cycle $\Gamma $ on $X\times Y$ we will denote by
$\Gamma ^t$ its transpose lying on $Y\times X$. By $M_{\sim
}:\mathcal V^{opp}_k\to \mathcal M_{\sim }$ we will denote the
functor which associates to any $X\in \mathcal V_k$ its motive
$M_{\sim }(X) =(X,id,0)$, where $id=[\Delta _X]$ is the class of
the diagonal $\Delta _X$ in $Corr^0_{\sim }(X,X)$, and to a
morphism $f:X\to Y$ the correspondence $M(f)=[\Gamma _f^t]\in
Corr^0_{\sim }(Y,X)=Hom_{\mathcal M_{\sim }}(M(Y),M(X))$, where
$\Gamma _f=\{ (x,f(x))\mid x\in X\} \subset X\times Y$ is the
graph of $f$.

In the following we fix a {\it Weil cohomology} theory with
$L$-coefficients $H^*$, where $L$ is a field of characteristic
zero, see \cite{Kl} for the definition. For example, if $k$ is an
arbitrary field one can take the \'etale cohomology groups
$H^*_{et}(\bar X,\mathbb Q_l)$ over the algebraic closure $\bar
k$, i.e. $\bar X=X\times _{k}\bar k$ and $l\ne char(k)$. If
$k=\mathbb C$ one can take also the usual Betti cohomology. Then
one defines a functor $H^i:\mathcal M_{rat}\to Vect_L$ for every
$i\in \mathbb Z$ by $H^i(M)= p_*H^{i+2m}(X)$ where $M=(X,p,m)$. By
$cl:A^i_{rat}(X)\to H^{2i}(X)$ we denote the cycle map; then
$\alpha \in A^i_{rat}(X)$ is homologically equivalent to zero iff
$cl(\alpha )=0$.

If $\sim $ is rational equivalence then $\mathcal M_{rat}$ is
called the category of {\it Chow motives} over $k$ with
coefficients in $\mathbb Q$. In the following we will write $A^i$
for $A^i_{rat}$, $h(-)$ for the functor $M_{hom}(-):\mathcal
V_k\to \mathcal M_{hom}$ and $A^i(X)_{hom}$ for the kernel of
$cl$, i.e. the subgroup in $A^i(X)$ of cycles which are
homologically trivial.

Under these assumptions one may consider the following equivalence
relations $\sim $ on cycles: ({\it rat}) rational equivalence;
({\it alg}) algebraic equivalence; ({\it hom}) homological
equivalence and ({\it num}) numerical equivalence \cite{J2}. It is
known that
   $$
   ({\it rat})\Rightarrow ({\it alg})\Rightarrow ({\it hom})
   \Rightarrow ({\it num})
   $$
Rational equivalence is strictly finer than algebraic equivalence
already for divisors on curves; a famous counterexample by
Griffiths showed that algebraic equivalence is strictly finer than
homological equivalence, even modulo torsion, for codimension $2$
cycles on a complex $3$-fold. According to Grothendieck's Standard
Conjectures on algebraic cycles \cite{Kl} homological equivalence
and numerical equivalence should coincide. By a result of Jannsen
\cite{J3} the category $\mathcal M_{\sim}$ is abelian semisimple
iff $\sim $ is the numerical equivalence.

Now let $X\in \mathcal V_k$ and assume, for simplicity, that $X$
is irreducible of dimension $d$. If we suppose that the conjecture
$C(X)$ holds, see \cite{Kl}, p.14, i.e. if the K\"unneth
components $\Delta (i,2d-i)$ of the diagonal $\Delta _X$ are
algebraic (which is known to be true for curves, surfaces and
abelian varieties), then the idempotent $H^*(X)\to H^i(X)\to
H^*(X)$ is represented by an algebraic correspondence $\sigma _i$
which is an idempotent in $A^d_{hom}(X\times X)$. Therefore we get
a natural decomposition:
   $$
   h(X)\simeq \bigoplus _{1\le i \le 2d}h^i(X)
   $$
where $h^i(X)=(X,\sigma _i,0)$ and $\sigma _i$ is, in fact, the
K\"unneth component $\Delta (i,2d-i)$.

Following \cite{Mu1} we will say that $X$ has a {\it
Chow-K\"unneth decomposition} if there exist orthogonal
idempotents $\pi _i$ ($0\le i\le 2d$) in $A^d(X\times X)$, such
that $cl(\pi _i)=\Delta (2d-i,i)$ and
   $$
   [\Delta _X]=\sum _{0\le i \le 2d}\pi _i
   $$
in $A^d(X\times X)$. This implies that in $\mathcal M_{rat}$ the
motive $M(X)$ decomposes as follows:
   $$
   M(X)=\bigoplus _{0\le i\le 2d}M^i(X)
   $$
where $M_i(X)=(X,\pi _i,0)$.

Murre conjectured, see \cite{Mu1} and \cite{Mu2}, that every $X\in
\mathcal V_k$ has a Chow-K\"unneth decomposition over the
algebraic closure $\bar k$. This conjecture is true for curves,
surfaces, abelian varietes, uniruled threefolds and elliptic
modular varieties, see \cite{DMu} and \cite{dAMSt} for further
references. If $X$ and $Y$ have a Chow-K\"unneth decompositon,
then the same holds for $X\times Y$.

If $X$ is a smooth projective variety of dimension $d$ satisfing
conjecture $C(X)$ then, with the above notations, one has the
following isomorphisms \cite{J1}:
   $$
   A^d_{hom}(X\times X)=\bigoplus _iEnd_{\mathcal M_{hom}}(h^i(X))
   $$
   $$
   A^d(X\times X)=\bigoplus _iEnd_{\mathcal M_{rat}}(M_i(X))
   $$
In order to relate rational equivalence with homological
equivalence for algebraic cycles it is therefore natural to ask
under which conditions the map
   $$
   End_{\mathcal M_{rat}}(M_i(X))\to End_{\mathcal M_{hom}}(h^i(X))
   $$
(induced by the functor $h$) is an isomorphism. This is in turn
strictly related, see \cite{J1}, Prop 5.8, with the existence of a
suitable filtration on the Chow ring of $X\times X $, such that
the associated graded groups only depend on the motives $h^i(X)$,
or, equivalently, to Murre's Conjecture (see Section
\ref{filtrations} below).

In this paper we show how finite dimensionality of the motive
$M(X)$ (see Def. \ref{Def2.1}) is related with the existence of
such a filtration.

The paper is organized as follows: in Section \ref{filtrations} we
recall the Conjectures of Beilinson and Murre on the existence of
a suitable filtration $F^{\bullet }$ on the Chow ring of a smooth
projective variety, and then we relate them with Bloch's Conjecure
for surfaces.

In Section \ref{S2}, after recalling the definitons and properties
of finite dimensional motives and some results of \cite{AK}, we
prove Theorems \ref{th5} and \ref{th7} which relate the finite
dimensionality of the motive with Murre's Conjecture.

In Section \ref{S3} we show that for a smooth projective surface
over an algebraically closed field of characteristic $0$ with
$p_g=0$ the motive $M(X)$ is finite dimensional iff the Chow group
of $0$-cycles of $X$ is finite dimensional in the sense of
Mumford.

{\sc Acknowledgements.} The authors wish to thank Jacob Murre and
Ivan Panin for many useful comments on an early version of this
paper.

\section{The Conjectures of Beilinson, Bloch \\ and Murre}
\label{filtrations}

Beilinson has conjectured the existence of decreasing filtrations
on Chow groups of all smooth projective varieties over $k$ which
are uniquely determined by the action of correspondences on
algebraic cycles \cite{J1}:

\begin{conjecture}
\label{beilinson}
For every $X\in \mathcal V_k$ there exists a
decreasing filtration $F^{\bullet}$ on $A^i(X)$, such that:

(a) $F^0A^j(X)=A^j(X)$; $F^1A^j(X)=(A^j(X))_{hom}$;

(b) $F^{\bullet}$ is compatible with the intersection product of
cycles;

(c) $F^{\bullet}$ is compatible with $f^*$ and $f_*$ if $f:X\to Y$
is a morphism;

(d) (if the K\"unneth components of $\Delta _X$ are algebraic) the
associated graded group $Gr^{\nu }_FA^j(X)=F^{\nu }A^j(X)/F^{\nu
+1}A^j(X)$ depends only on the motive $h^{2j-\nu }(X)$ of $X$ in
$\mathcal M_{hom}$;

(e) $F^{j+1}A^j(X)=0$ for all $j$.
\end{conjecture}

If such a filtration exists then it is unique \cite{J1}. If the
K\"unneth components of the diagonal are algebraic, then a weaker
form (see \cite{J2}, p. 12) of the condition (d) is:

(${\rm d}^{\prime }$) Let $Y\in \mathcal V$ and let $\gamma \in
Corr^{j-i}(Y\times X)$. If the induced map $\gamma _*$ between
$H^{2i-\nu }(Y)$ and $H^{2j-\nu }(X)$ is zero, then so is the map
$Gr^{\nu }_F\gamma :Gr^{\nu }_FA^i(Y)\to Gr^{\nu }_FA^j(X)$.

Conjecture \ref{beilinson} is in turn equivalent (again assuming
that the K\"unneth components of the diagonal are algebraic) to
the following Conjecture of Murre, see \cite{Mu1} and \cite{J1}:

\begin{conjecture}
\label{murre} For any smooth projective (irreducible, for
simplicity) variety $X$ of dimension $d$:

(I) there exists a Chow-K\"unneth decomposition $[\Delta _X]=\sum
_{i=0}^{2d}\pi _i$;

(II) the correspondences $\pi _0,\dots ,\pi _{j-1}$ and $\pi
_{2j+1},\dots ,\pi _{2d}$ act as $0$ on $A^j(X)$;

(III) if $F^{\bullet }$ is the filtration on $A^j(X)$ defined by
$F^{\nu }A^j(X)=ker(\pi _{2j})\cap ker(\pi _{2j-1})\cap \dots \cap
ker(\pi _{2j-\nu +1})$, then $F^{\bullet }$ is independent of the
choice of the projectors $\pi _i$;

(IV) $F^1A^j(X)=A^j(X)_{hom}$.
\end{conjecture}

The status of Conjecture \ref{murre} is as follows: it is
trivially true for curves. For surfaces and for the product of a
surface with a curve Murre has shown the existence of a
Chow-K\"unneth decomposition satsfying (II) and (IV), see
\cite{Mu3}, \cite{Mu1} and \cite{Mu2}. For surfaces he also shows
that there is a filtration which is the natural one, i.e. it
coincides with the filtration for $0$-cycles considered in
\cite{B}. For abelian varieties the existence of a Chow-K\"unneth
decomposition follows from works of Shermenev, Denninger-Murre and
K\"uenneman (see \cite{Ku} for references): part of (II) is true
(see Remark \ref{abvar}) and if (II) is true then (III) is also
true for a natural choice of the projectors $\pi _i$ \cite{Mu1}.

Let us consider the case when $X$ is a smooth projective surface
over an algebraically closed field $k$. By results in \cite{Mu3},
$X$ has a Chow-K\"unneth decomposition $[\Delta _X]=\sum
_{i=0}^4\pi _i$, where $\pi _0=[x_0\times X]$ and $\pi _4=[X\times
x_0]$ are the trivial projectors induced by a fixed point $x_0\in
X$, $\pi _1$ is the Picard projector (which is closely connected
with the Picard variety $Pic^0(X)$ of the surface $X$), $\pi
_3=\pi _1^t-\pi _1\circ \pi _1^t$ is the Albanese projector
(connected with the Albanese variety $Alb(X)$ of $X$) and $\pi
_2=\Delta _X -\pi _0-\pi _1-\pi _3-\pi _4$. The projectors $\pi
_i$ yield the motivic decomposition
  $$
  M(X)=\sum _{0\le i\le 4}M_i(X)\; ,
  $$
where $M_i(X)=(X,\pi _i,0)$ for any $i=\{ 0,\dots ,4\}$, and the
corresponding Murre's filtration is:
     $$
     F^0A^i(X)=A^i(X)\; ,
     $$
     $$
     F^{i+1}A^i(X)=0\; ,
     $$
     $$
     F^1A^1(X)=A^1(X)_{hom}=A^1(X)_{num}=\ker (\pi _2)\; ,
     $$
     $$
     F^1A^2(X)=ker(\pi _4)=A^2(X)_0
     $$
-- the group of zero-cycles of degree $0$ on $X$, and
     $$
     F^2A^2(X)=\ker (\pi _3\! \mid _{F^1})=T(X)\; ,
     $$
where $T(X)$ is so called Albanese kernel of the surface $X$, i.e.
the kernel of the Abel-Jacobi map $a_X:A^2(X)_0\to Alb(X)$. The
graded group $Gr^*_F(A^2(X))$ associated to the filtration above
is:
   $$
   \mathbb Q\oplus Alb(X)_{\mathbb Q}\oplus T(X)\; .
   $$

A similar (truncated) filtration  $F^{\bullet}$ can be defined for
the Chow group of $0$-cycles of any smooth variety $Y$ of
dimension $d$. Then one has, in analogy to the case of surfaces:
   $$
   Gr^*_FA^d(Y)=\mathbb Q\oplus Alb(Y)_{\mathbb Q}\oplus T(Y)\; .
   $$

If Beilinson's Conjectural Filtration $F^{\bullet }$ exists for
every smooth projective variety and $X$ is a surface, then any
correspondence $\gamma \in A^2(Y\times X)$, where $d=dim(Y)$,
respects the filtration; if $\gamma \in A^2(Y\times X)_{hom}$ then
$\gamma _*$ acts as $0$ on $Gr^*_FA^d(Y)$. This shows that the
Beilinson's Conjecture implies the following conjecture formulated
in \cite{B}:

\begin{conjecture}
\label{bloch1} Let $X$ be a smooth projective surface and let $Y$
be a smooth projective variety of dimension $d$. For any $\gamma
\in A^2(Y\times X)$ its action on $Gr^*_FA^d(Y)$:
   $$
   Gr^*_F\gamma :Gr^*_FA^d(Y)\longrightarrow Gr^*_FA^2(X)
   $$
depends only upon the cohomology class $cl(\gamma )$ in
$H^4(Y\times X)$.
\end{conjecture}

Conjecture \ref{bloch1} implies

\begin{conjecture}
\label{bloch2} If $X$ is a complex surface with geometric genus
$p_g=0$, then the Albanese kernel $T(X)$ vanishes, see \cite{B},
1.11.
\end{conjecture}

Note that, by a result of \cite{Ro}, if $k$ is algebraically
closed then the kernel of the Abel-Jacobi map $a_X$, considering
with coefficients in $\mathbb Z$, is torsion free.

Bloch's conjecture on  the Albanese kernel holds for surfaces of
Kodaira dimension less than $2$ \cite{BKL} and it is still open
for complex surfaces of general type with $p_g=0$, see
\cite{InMiz}, \cite{Voi} and \cite{GP1}.

\begin{remark}
{\rm In general Bloch's Conjecture \ref{bloch1} does not imply
that the action $\gamma :A^d(Y)\to A^2(X)$ only depends on the
cohomology class $cl(\gamma )$. In fact, let $X$ be a complex
surface with $p_g(X)>0$ and $q(X)=0$ (where
$q(X)=dim(H^1(X,\mathcal O_X))$ is the  irregularity of $X$) and
let $C$ be a generic curve on $X$. Let $Y=S^2C$ be the symmetric
square of the curve $C$. Then $A^d(Y)\simeq J(C)\oplus T(Y)$ and
$A^2(X)\simeq T(X)$, where $J(C)$ is the Jacobian of the curve
$C$, $T(Y)$ and $T(X)$ are the Albanese kernels. The map
$f:S^2C\to S^2X$ yields a series of effective $0$-cycles of degree
$2$ on $X$. Let $\Gamma \subset Y\times X$ be the associated
correspondence, i.e.
   $$
   [\Gamma ]=\{ [(Y,P+Q)]\mid P,Q\in C\}
   \subset A^2(Y \times X)
   $$
and let $\gamma =[\Gamma ]$. Then the class $cl(\gamma )$ in
$H^4(Y\times X)$ has components $\gamma(0,4)$, $\gamma (4,0)$ and
$\gamma(2,2)$. By adding constant correspondences to $\gamma $ we
may assume that $\gamma (0,4)= \gamma (4,0)=0$. Moreover the
component $\gamma (2,2)$ in $H^2(Y)\otimes H^2(X)$ belongs to
$NS(Y)\otimes NS(X)$. Therefore the action of $\gamma (2,2)$ on
$0$-cycles is trivial because every $0$-cycle can be moved away
from a finite number of divisors. The graded map
   $$
   Gr^*_F\gamma :Gr^*_FA^d(Y)\to Gr^*_FA^2(X)
   $$
is $0$. In fact we have $\Gamma \subset Y\times C\subset Y\times
X$, whence $\gamma $ can also be viewed as a correspondence
between $Y$ and $C$. As such it determines a map
   $$
   \gamma ^{\prime }:A^d(Y)\simeq J(C)\oplus T(Y)\to J(C)\; ,
   $$
which is just the projection onto the first factor. Since $\gamma
$ factors trough $\gamma ^{\prime }$, we see that $\gamma $ is $0$
on $T(Y)$ and, therefore, $Gr^*_F\gamma $ is the zero map. However
the map
   $$
   \gamma :A^d(Y)\to A^2(X)=T(X)
   $$
is not zero: in fact, $C$ being a general curve on the surface $X$
with $p_g(X)>0$, the map induced by $\gamma $ between $J(C)$ and
$T(X)$ is non trivial. This is the consequence of a  famous
results of Mumford on the group of $0$-cycles on surfaces with
$p_g>0$, see \cite{Voi}, pg. 186.}
\end{remark}

\section{Finite dimensional motives and \\ Murre's Conjecture}
\label{S2}

In this section we first recall the definition and some results on
finite dimensional motives, which have been introduced by S.-I.
Kimura in \cite{Kim}, and then prove our results relating finite
dimensionality with the Conjectures stated in Section
\ref{filtrations}.

Let $\mathcal C$ be a pseudoabelian, $\mathbb Q$-linear, tensor
category and let $X$ be an object in $\mathcal C$. Let $\Sigma _n$
be the symmetric group of order $n$. Any $\sigma \in \Sigma_n$
defines an endomorphism $\Gamma _{\sigma }:(x_1,\dots ,x_n)\mapsto
(x_{\sigma (1)},\dots ,x_{\sigma (n)})$ of the $n$-fold tensor
product $X^n$ of $X$ by itself. There is a one-to-one
correspondence between all irreducible representations of the
group $\Sigma _n$ (over $\mathbb Q $) and all partitions of the
integer $n$. Let $V_{\lambda }$ be the irreducible representation
corresponding to a partition $\lambda $ of $n$ and let $\chi
_{\lambda }$ be the character of the representation $V_{\lambda
}$. Let
   $$
   d_{\lambda }=\frac{\dim (V_{\lambda })}{n!}\sum_{\sigma \in
   \Sigma _n}\chi _{\lambda }(\sigma )\cdot \Gamma _{\sigma }
   $$
Then $\{ d_{\lambda }\} $ is a set of pairwise orthogonal
idempotents in $End_{\mathcal C}(X^n)$, such that $\sum d_{\lambda
}=id_{X^n}$. The category $\mathcal C$ being pseudoabelian they
give a decomposition of $X^n$. The $n$-th symmetric product $S^nX$
of $X$ is then defined to be $im(d_{\lambda })$ when $\lambda $
corresponds to the partition $(n)$, and the $n$-th exterior power
$\wedge ^nX$ is $im(d_{\lambda })$ when $\lambda $ corresponds to
the partition $(1,\dots ,1)$. In particular, we have symmetric and
exterior powers in $\mathcal M_{\sim }$.

The following definition was made in \cite{Kim}, see also
\cite{GP1} or \cite{AK}:

\begin{definition}
\label{Def2.1} The object $X$ in $\mathcal C$ is said to be {\it
evenly (oddly) finite dimensional} if $\wedge ^nX=0$ ($S^nX=0$)
for some $n$. An object $X$ is finite dimensional if it can be
decomposed into a direct sum $X_+\oplus X_-$ where $X_+$ is evenly
finite dimensional and $X_-$ is oddly finite dimensional.
\end{definition}

Now we want to show that, if the motive $M(X)$ is finite
dimensional, then $X$ has a Chow-K\"unneth decomposition. We first
recall a result which has been proved in \cite{J1}, 5.3:

\begin{lemma}
\label{lemma1} Assume $X$ is a smooth projective variety of
dimension $d$, such that $A^d(X\times X)_{hom}$ is a nilpotent
ideal of $A^d(X\times X)$. Assume moreover that the K\"unneth
components of the diagonal are algebraic. Then $X$ has a
Chow-K\"unneth decomposition.
\end{lemma}

\begin{theorem}
\label{th2} Let $M$ be a finite dimensional motive in $\mathcal
M_{rat}$ and let $f$ be a homologically trivial endomorphism of
$M$, i.e. $f$ induces the $0$ map on $H^*(M)$. Then $f$ is
nilpotent in $End_{\mathcal M_{rat}}(M)$.
\end{theorem}

\begin{pf}
See \cite{Kim}, 7.2
\end{pf}

\begin{corollary}
\label{corollary3} Let $M(X)$ be a finite dimensional Chow motive.
Assume that the K\"unneth components of the diagonal of $X$ are
algebraic. Then $X$ has a Chow-K\"unneth decomposition.
\end{corollary}

\begin{pf}
Apply Theorem \ref{th2} and Lemma \ref{lemma1}
\end{pf}

\begin{remark}
{\rm If $M(X)$ has a Chow-K\"unneth decomposition then the
projectors $\pi _i$ defining the motives $M_i(X)$ are by no means
unique: for instance the cycle class of the trivial projector $\pi
_0$ depends on the choice of a rational point $x_0$ on $X$.
However the motives $M_0(X)$ and $M_{2d}$ are unique (up
isomorphisms in $\mathcal M$). Also, for a curve $C$, uniqueness
of the motives $M_i(X)$ for $i=0,1,2$ is easy \cite{Mu3}, 5.1. For
an arbitrary $X$ of dimension $d$ Murre has shown \cite{Mu3}, 5.2,
that the motives $M_1(X)=(X,\pi _1,0)$ and $M_{2d-1}(X)=(X,\pi
_{2d -1},0)$, where $\pi _1$ and $\pi _{2d-1}$ are respectively
the Picard and the Albanese projectors, are, up to isomorphisms,
independent of the polarization choosen to construct $\pi _1$ and
$\pi _{2d -1}$.

We will show in Theorem \ref{th5} that, if $M(X)$ is finite
dimensional, then all the $M_i(X)$ are unique, up to
isomorphisms.}
\end{remark}

The main known properties of finite dimensional objects are:

1) If two objects $X,Y\in \mathcal C$ are finite dimensional so is
their direct sum $X\oplus Y$ and their tensor product $X\otimes
Y$. If $X$ is a subobject of a finite dimensional object $Y$ then
$X$ is finite dimensional (equivalently, if $X$ is a quotient
object of a finite dimensional object $Y$, it is finite
dimensional). Moreover, a direct summand of an evenly (oddly)
finite dimensional motive is evenly (oddly) finite dimensional.
Note that these properties were proved by Kimura for Chow motives
over a field. But they can be proved in an arbitrary pseudoabelian
$\mathbb Q$-linear tensor category,see \cite{AK}.

2)In particular the properties in 1) impliy the following .If $f:Y\to X$ is a
proper surjective morphism of smooth projective varieties and
$M_{\sim }(Y)$ is finite dimensional then $M_{\sim }(X)$ is also
finite dimensional; the motoive $M_{\sim }(X)\otimes M_{\sim
}(Y)=M_{\sim }(X\times Y)$ of the fibered product $X\times Y$ is
finite dimensional if $M_{\sim }(X)$ and $M_{\sim }(Y)$ are finite
dimensional.

3) If a motive $M$ is evenly and oddly finite dimensional then
$M=0$ \cite{Kim}, 6.2.

4) The dual object $X^*$ in a rigid category $\mathcal C$ is
finite dimensional iff $X$ is finite dimensional.

5) Finite dimensionality is a birational invariant for surfaces,
\cite{GP1}, Th. 2.8.

The following theorem gives classes of smooth projective varieties
whose motives are finite dimensional

\begin{theorem}
\label{prop4} (i) The motive of a smooth projective curve over a
field is finite dimensional. (ii) The motive of a variety which is
the quotient of a product $C_1\times\cdots C_n$ of curves under
the action of a finite group $G$ acting freely on $C_1\times\cdots
\times C_n$ is finite dimensional. (iii) If $X$ is an abelian
variety, then $M(X)$ is finite dimensional. (iv) The same result
holds if $X$ is a Fermat hypersurface of degree $d$ in $\mathbb
P^n$.
\end{theorem}

\begin{pf}
(i) was proved in \cite{Kim}. (ii) and (iii) follow from (i) and
the above properties 1) -- 5). For abelian varieties see also
\cite{Sch}, 3.4. The proof of the fact that the motive of a Fermat
hypersurface is finite dimensional can be found in \cite{GP1}.
\end{pf}

Let $\mathcal M_{Kim}$ be the full subcategory of $\mathcal
M_{rat}$ generated by finite dimensional oblects. From the
properties 1) -- 5) it follows then that $\mathcal M_{Kim}$ is a
pseudoabelian, rigid and tensor category. Kimura stated

\begin{conjecture}
$\mathcal M_{Kim}=\mathcal M_{rat}$
\end{conjecture}

Evidently, $\mathcal M_{Kim}$ contains a subcategory generated by
the Chow motives of varieties  as in  Theorem \ref{prop4}, their products
and quotients in $\mathcal M_{rat}$.

The relations between finite dimensionality and the Conjectures
stated in Section \ref{filtrations} can be made more precise using
some recent results from \cite{AK}. We first recall the definition
of {\it dimension} for an object in a rigid tensor category
$\mathcal C$, see \cite{AK} or \cite{DeMi}.

For any $X\in \mathcal C$ let $\epsilon_X:X\otimes X^*\to 1$ be
the evaluation map, and for any two $X,Y\in \mathcal C$ let
   $$
   i_{X,Y}:Hom_{\mathcal C}(1,X^*\otimes Y)
   \stackrel{\simeq }{\longrightarrow }
   Hom_{\mathcal C}(X,Y)
   $$
be the canonical isomorphism. Let $h\in End_{\mathcal C}(X)$. Then
we define the trace of $h$ to be
   $$
   tr(h)=\epsilon _{X^*}\circ i^{-1}_{X,X}(h)
   \in Hom_{\mathcal C}(1,1)\simeq \mathbb Q
   $$
and define $dim(X)=tr(id_X)$.

If $dim(X)=d$ then
   $$
   dim(\wedge ^nM)={d\choose n}={d(d-1)\cdot \dots \cdot (d-n+1)\over n!}
   $$
and
   $$
   dim(S^nA)=
   {d+n-1\choose n}={d(d+1)\cdot \dots \cdot (d+n-1\over n!}\; ,
   $$
see \cite{AK}, 7.2.4. Therefore, if $dim(X)=d>0$, then $dim(\wedge
^{d+1}X)=0$; if $dim(X)=-d<0$ then $dim(S^{d+1}X)=0$.

This dimension  is related  with Kimura's finite dimensionality
in the following way.

\begin{definition}
\label{kim} Let  $X\in \mathcal C$ be afinite dimensional object.
Then  $kim(X)$ is   the smallest
integer $n$, such that $\wedge^nX=0$ if $X$ is evenly finite
dimensional, and   $S^nX=0$ if $X$ is oddly finite
dimensional.
\end{definition}

If $X$ is Kimura finite dimensional, then $dim(X)$ is an integer
\cite{AK}, 9.1.5: if $X$ is evenly finite dimensional then
$dim(X)=kim(X)$, while if $X$ is oddly finite dimensional then
$dim(X)=-kim(X)$.

If $H$ is a Weil cohomology theory (with coefficients in a field
$L$ of characteristic zero) on $\mathcal V_k$ then for every Chow
motive $M\in \mathcal
M_{rat}$ we have $dim(M)=\sum _{i\in \mathbb Z} (-1)^idim(H^i(M))$.
For all $X\in \mathcal V_k$
which satisfy the standard conjecture $C(X)$, i.e. the K\"unneth components of
the diagonal $\Delta _X$ are algebraic, there exist projectors
$p^+_M $ and $p^-_M$ in $End_{\mathcal M_{rat}}(M(X))$, such that
$H(p^+_M)$ and $H(p^-_M)$ are the projectors corresponding to the
splitting of $H(X)$ respectively into the even and the odd part.

Let $\mathcal A$ be the full subcategory of $\mathcal M_{rat}$ of
objects $A$, such that projectors $p^+_A$ and $p^-_A$ exist in
$End_{\mathcal M_{rat}}(A)$. Then $\mathcal A$ is a rigid, tensor
and $\mathbb Q$-linear subcategory of $\mathcal M_{rat}$
containing all the motives of curves, surfaces, abelian varieties,
their products and subobjects. For every object $A\in \mathcal A$
the projectors $p^+_A$ and $p^-_A$ induce a decomposition
$A=A^+\oplus A^-$, see \cite{AK}, 8.3.

If $A\in \mathcal A$ has a decomposition $A=A^+\oplus A^-$ then
$A^+$ is evenly finite dimensional and $A^-$ is oddly finite
dimensional (and hence $A$ is finite dimensional) iff there exists
an integer $n$, such that :
   $$
   s\wedge ^nA=0\; ,
   $$
where $s\wedge ^nA=\bigoplus _{i+j=n}\wedge ^iA^+\otimes S^jA^-$.
If such $n$ exists then the smallest one is the integer
$kim(A^+)+kim(A^-)+1$, see \cite{AK}, 9.1.11.

If $A$ is finite dimensional then the decomposition $A=A^+\oplus
A^-$ is unique up to isomorphisms, i.e. if $A=\tilde A^+\oplus
\tilde A^-$, where $\tilde A^+$ ($\tilde A^-$) is evenly (oddly)
finite dimensional, then $A^+\simeq \tilde A^+$ and $A^-\simeq
\tilde A^-$, see \cite{Kim}, 6.3.

\begin{theorem}
\label{th5} Let $X$ be a smooth projective variety over $k$, such
that the K\"unneth components of the diagonal $\Delta _X$ are
algebraic. Assume that the motive $M(X)\in \mathcal M_{rat}$ is
finite dimensional. Then $M(X)$ has a Chow-K\"unneth decomposition
   $$
   M(X)=\bigoplus _{0\le i\le 2d}M_i(X)
   $$
with $M_i(X)=(X,\pi _i,0)$, which is independent of the choice of
the projectors $\pi _i$, i.e., if $\{ \tilde \pi _i\} $ is another
set of orthogonal idempotents lifting the K\"unneth components of
$\Delta _X$, then
   $$
   M_i(X)\simeq \tilde M_i(X)
   $$
in $\mathcal M_{rat}$, where $\tilde M_i(X)=(X,\tilde \pi
_i,0)$.
\end{theorem}

\begin{pf}
By Corollary \ref{corollary3} the motive $M(X)$ has a
Chow-K\"unneth decomposition. Let $M(X)=\bigoplus _{0\le i\le 2d}
M_i(X)$, where $M_i(X)=(X,\pi _i,0)$ and let $\{ \tilde \pi _i\} $
be another set of orthogonal idempotents lifting the K\"unneth
components of $\Delta _X$.

Let's consider the following composition of projectors, for $i=
0,\dots ,2d$:
   $$
   M_i(X)
   \stackrel{\pi _i}{\rightarrow}
   M(X)
   \stackrel{\tilde \pi_i}{\rightarrow}
   M_i(X)
   \stackrel{\tilde \pi_i}{\rightarrow}
   M(X)
   \stackrel{\pi _i}{\rightarrow}
   M_i(X)
   $$
and set
   $$
   e_i=\pi _i\circ \tilde \pi _i\circ \tilde \pi _i\circ \pi _i=
   \pi _i\circ \tilde \pi _i\circ \pi _i\; .
   $$
Then $e_i\circ \pi _i=\pi _i\circ e_i$, i.e. $e_i\in End_{\mathcal
M_{rat}}(M_i(X))$.

We claim that $e_i= \pi _i$, i.e. $e_i$ is the identity on
$M_i(X)$.

$M(X)$ being finite dimensional from Th. \ref{th2} it follows that
$I=A^2(X\times X )_{hom}$ is a nilpotent ideal of $A^2(X\times
X)$. Therefore there exists an element $\eta \in I$, such that
$\tilde \pi _i=(1+\eta )^{-1}\circ \pi _i\circ (1+\eta )$ for $i
=0,\dots ,2d$, see \cite{J1}, 5.4, and we have:
    $$
    e_i-\pi _i=\pi _i-(\pi _i-\pi _i\circ \eta \circ \pi _i\circ \eta
    \circ \pi _i)=\pi _i\circ \eta \circ \pi _i\circ \eta \circ \pi _i\; .
    $$
So we are left to show that $\pi _i\circ \eta \circ \pi _i\circ
\eta \circ \pi _i=0$. By induction on the index of nilpotency of
$I$ we may assume that $I^2=0$. Then we can take $\eta =\pi _i
\circ \epsilon _i-\epsilon _i\circ \pi _i$ where $\tilde \pi _i=
\pi _i+\epsilon _i$ with $\epsilon_i \in I$ and $\epsilon _i^2=0$,
see \cite{Mu3}, page 203. Expanding $(\pi _i+\epsilon _i)^2$ leads
to the equation $\epsilon _i=\pi _i\circ \epsilon _i +\epsilon
_i\circ \pi _i$, whence:
   $$
   \pi _i\circ \epsilon _i\circ \pi _i=\epsilon _i\circ \pi _i\epsilon _i=0\; .
   $$
 From the equalities above we get:
   $$
   \pi _i\circ \eta \circ \pi _i\circ \eta \circ \pi _i=
   \pi _i\circ (\pi _i\circ \epsilon _i-\epsilon _i\circ \pi _i)
   \circ \pi _i\circ (\pi _i\circ \epsilon _i-\epsilon _i\circ
   \pi _i)\circ \pi _i=
   $$
   $$
   \pi _i\circ \epsilon _i\circ \pi _i\circ
   \epsilon _i\circ \pi _i=0\; .
   $$
In a completely similar way one shows that $\tilde e_i=\tilde \pi
_i\circ \pi _i\circ \tilde \pi _i$ is the identity on $\tilde
M_i(X)$. Therefore, $\tilde \pi _i\circ \pi _i$ yields an
isomorphism $M_i(X)\simeq \tilde M_i(X)$.
\end{pf}

\begin{remark}
\label{abvar} (Abelian varieties) {\rm Let $X$ be an abelian
variety of dimension $d$ over an algebraically closed field $k$ of
char $0$. Then $M(X)$ has a Chow-K\"unneth decomposition;
moreover, there exists a unique decompositon $[\Delta_X ]= \sum
_i\pi _i\in A^d(X \times X)$, such that
   $$
   n^*\circ \pi _i=\pi _i\circ n^*=n^i\pi _i
   $$
for every $n\in \mathbb Z$, where $n^*=(id_X\times n)^*$ and $n$
is the multiplication by $n$ on $X$. The correspondenses $\{ \pi
_i\} $ are orthogonal projectors, such that $\pi _0,\dots , \pi
_{j-1}$ and $\pi _{j+d+1},\dots ,\pi _{2d}$ operates as $0$ on
$A^j(X)$, see \cite{Mu1}, 2.5.2. The corresponding decomposition
$M(X)=\sum _{0\le i\le 2d} M_i(X)$ satisfies a part of conditon
(II) in Conjecture \ref{murre}. The motive $M(X)$ is finite
dimensional: from Theorem \ref{th5} it follows that this
decomposition is unique (up to isomorphism). Therefore, if there
exists a Chow-K\"unneth decompositon satisfing the rest of the
condition (II), i.e. such that also $\pi _i$ operates as $0$ on
$A^j(X)$ for $2j+1\le i\le j+d$, then it is isomorphic to the one
above. This condition is in turn equivalent to {\it Beauville's
Conjecture}, see \cite{Mu1}, 2.5.3, and \cite{Bea}, on the
vanishing of the groups $A^j_s(X)=\{ \alpha \in A^j(X)\mid
n^*\alpha =n^{2j-s}\alpha \} $ for $s<0$.

Beauville's Conjecture being true for cycles of codimension
$j=0,1,d-2,d-1$ it follows that conditon (II) is in particular
satisfied for all abelian varieties of dimension at most $4$.
Therefore, for all abelian varieties which satisfy Beauville's
Conjecture, the filtration associated to a Chow -K\"unneth
decomposition is independent of the choices of the projectors, in
the sense that it only depends on the isomorphism classes of the
motives $M_i(X)$. This proves that Beauville's Conjecture implies
Murre's conjecture for an abelian variety.}
\end{remark}

In \cite{AK}, 9.2.4, it has been remarked that if Beilinson's
Conjecture or, equivalently, Conjecture \ref{murre} is true for
all varieties $X$ and also the Standard Conjectures hold, then all
Chow motives of smooth projective varieties are finite
dimensional, i.e. Kimura's Conjecture holds. The following Theorem
\ref{th7} avoids the assumption about the Standard Conjectures.

We first prove a lemma which is a direct consequence of a result
in \cite{J1}, 5.8.

\begin{lemma}
\label{lemma6} Let $Y$ be a smooth projective variety of dimension
$d$, such that $Y$ has a Chow-K\"unneth decomposition, say
$[\Delta _Y]=\sum _{i=0}^{2d}\pi _i$, and $Y\times Y$ satisfies
the Murre Conjecture. Let  $N_i=(Y,\pi _i,0)$: then
    $$
    Hom_{\mathcal M_{rat}}(N_s,N_t)=0\; \; \hbox{if}\; \; s\ne t
    $$
and
    $$
    Hom_{\mathcal M_{rat}}(N_s,N_s)=Hom_{\mathcal M_{hom}}
    (h(N_s),h(N_s))
    $$
for any $s\in \{ 0,1,\dots ,2d\} $.
\end{lemma}

\begin{pf} Let $\tilde \pi _i=\pi _{2d -i}^t$ be the transpose of
$\pi _{2d -i}$ and let $\Pi _r= \sum _{i+j=r}\tilde \pi _i\times
\pi _j$. By the same argument as in \cite{J1}, 5.8, the projector
$\Pi _r$ is a lifting of the $r$th K\"unneth component of the
diagonal $\Delta _{Y\times Y}(4d-r,r)$. Since $Y\times Y$
satisfies the conditon (II) in Conjecture \ref{murre}, it follows
that $\Pi _r$ acts as $0$ on $A^d(Y\times Y)$ for $r>2d$, whence
we get, for all pairs $(i,j)$ with $i+j>2d$:
   $$
   0=(\tilde \pi _i\times \pi _j)A^d(Y\times Y)=
   \pi _j\circ Corr^0(Y,Y)\circ \tilde \pi _i^t=
   $$
   $$
   =\pi _j\circ Corr^0(Y,Y)\circ \pi _{2d-i}=
   Hom_{\mathcal M_{rat}}(N_{2d-i},N_j)\; .
   $$
This shows that $Hom_{\mathcal M_{rat}}(N_s,N_t)=0$ for $s<t$.

If we take $\tilde \pi _i=\pi _i^t$ and $\tilde \pi _j=\pi
_{2d-j}$, the projector $\Pi _r=\sum _{i+j=r}\tilde \pi _i\times
\tilde \pi _j$ is (up to an isomorphism of $H^*(Y\times Y \times Y
\times Y)$) again a lifting of $\Delta _{Y\times Y}(4d-r,r)$. As
such $\Pi _r$ acts as $0$ on $A^d(Y\times Y)$ for $r>2d$. Just as
before we get, for all pairs $(i,j)$ with $i+j>2d$:
    $$
    0=(\tilde \pi _i\times \tilde \pi _j)A^d(Y\times Y)=
    \pi _{2d-j}\circ Corr^0(Y,Y)\tilde \pi _i^t=
    $$
    $$
    =\pi _{2d-j}\circ Corr^0(Y,Y)\circ \pi _i=
    Hom_{\mathcal M_{rat}}(N_i,N_{2d-j})\; .
    $$
Therefore, $Hom_{\mathcal M_{rat}}(N_s,N_t)=0$ for $s>t$.

The proof of the equality $Hom_{\mathcal M_{rat}}(N_s,N_s)=
Hom_{\mathcal M_{hom}}(h(N_s),h(N_s))$ follows from the same
argument as in \cite{J1}, 5.8: one takes projectors $\Pi _r=\sum
_{i+j=r}\tilde \pi _i\times \pi _j$ where $\tilde \pi _i=\pi
_{2d-i}^t$ and applies conditon (IV) in Murre's Conjecture. Then
$A^d(Y\times Y)_{hom}=F^1A^d(Y\times Y)=ker(\Pi _{2d})$ and we
obtain
   $$
   (\tilde \pi _{2d-s}\times \pi _s)A^d(Y\times Y) =
   (\tilde \pi _{2d-s}\times \pi _s)Corr^0_{hom}(Y,Y)=
   $$
   $$
   =\Delta _Y(2d-s,s)\circ Corr^0_{hom}(Y,Y)\circ \Delta _Y(2d-s,s)=
   $$
   $$
   =Hom_{\mathcal M_{hom}}(h(N_s),h(N_s))\; .
   $$
This proves that $Hom_{\mathcal M_{rat}}(N_s,N_s)= Hom_{\mathcal
M_{hom}}(h(N_s),h(N_s))$.
\end{pf}

\begin{theorem}
\label{th7} Let $X$ be a smooth projective variety of dimension
$d$ over $k$. Let $n=\sum _idim(H^i(X))$ and let $m=n+1$. Assume
that $X$ has a Chow-K\"unneth decomposition and $X^m\times X^m$
satisfies Murre's Conjecture. Then the motive $M(X)$ is finite
dimensional.
\end{theorem}

\begin{pf} There exist projectors $p_+$ and $p_-$ splitting the
motive $M=M(X)$ into $M^+$ and $M^-$, such that the cohomology of
$M^+$ is $H^+(X)=\sum _{i\in \mathbb Z}H^{2i}(X)$ and the
cohomology of $M^-$ is $H^-(X)=\sum _{i\in \mathbb Z}H^{2i+1}(X)$.
Therefore, $M(X)$ is finite dimensional iff $M^+$ is evenly finite
dimensional and $M^-$ is oddly finite dimensional. We have:
$dim(M^+)=B_+=dim(H^+(X))$ and $dim(M^-)=-B_-=-dim(H^-(X))$.
Therefore,
   $$
   dim(s\wedge ^mM)=dim\left(\sum _{i+j=m}\wedge ^iM\otimes
   S^jM\right)=0
   $$
if $m=B_++B_-+1$. So, in order to show that $M$ is finite
dimensional, it is enough to prove that
   $$
   s\wedge ^mM=\bigoplus _{i+j=m}\wedge ^iM^+\otimes S^jM_-=0\; .
   $$

The functor $H:\mathcal M_{hom}\to Vect_L$ being faithful, from
$dim(s\wedge ^mM)=0$ we get: $s\wedge ^mh(M)=0$, see \cite{AK},
8.3.1. Let $q^+_i$ and $q^-_j$ be the projectors which define
respectively $\wedge ^iM^+$ and $S^iM^-$. Then the projector
$q=\sum _{i+j=m} q^+_i\otimes q^-_j$, which defines $s\wedge ^mM$,
belongs to $End_{\mathcal M_{rat}}(M(X^m))$ and is homologically
trivial, i.e. $h(q)=0$.

We claim that $q=0$, i.e. $s\wedge ^mM=0$.

Let $Y= X^m$. Since $X$ has a Chow-K\"unneth decomposition also $Y
=X^m$ has a Chow-K\"unneth decomposition, see \cite{Mu2}, 5.1.
Moreover, $Y\times Y$ satisfies Murre's Conjecture by assumptions.
Let $M(Y)=\sum _{0\le s\le 2md}N_i$ where $N_i=(Y,\pi _i,0)$ be a
Chow-K\"unneth decomposition for $Y$. From Lemma \ref{lemma6} it
follows that:
    \begin{equation}
    \label{tag7}
    Hom_{\mathcal M_{rat}}(N_s,N_t)=
    \left\{
    \begin{array}{ll}
    0 & \mbox{if $s\neq t$,} \\
    Hom_{\mathcal M_{hom}}(h(N_s),h(N_s)) &
    \mbox{if $s=t$}
    \end{array}
    \right.
    \end{equation}
Let $f_{s,t}=\pi _t\circ q\circ \pi _s\in Hom_{\mathcal
M_{rat}}(N_s,N_s)$ be the composition map:
    $$
    N_s
    \stackrel{\pi _s}{\rightarrow}
    M(Y)
    \stackrel{q}{\rightarrow}
    M(Y)
    \stackrel{\pi _t}{\rightarrow}
    N_t\; .
    $$
Then $\sum _s\pi _i\circ q=\sum _sq\circ \pi _s=q$ and $\sum
_{s,t}f_{s,t}=\sum _t\pi _t\circ \sum _sq\circ\pi _s=\sum _t\pi _t
\circ q=q$. Therefore we get:
   $$
   q=\sum _{s\ne t}(\pi _t\circ q\circ \pi _s)+\sum _s\pi _s\circ
   q\circ \pi _s\; .
   $$
 From (\ref{tag7}) it follows that $\sum _{s\ne t}(\pi _t\circ
q\circ \pi _s)=0$ which yields:
   $$
   q=\sum _{0\le s\le 2md}\pi _s\circ q\circ \pi _s\in
   Hom_{\mathcal M_{rat}}(N_s,N_s)
   $$
with $h(q)=0$. From the second equality in (\ref{tag7}) it follows
that $q=0$. This proves that $s\wedge ^mM=0$.
\end{pf}

\begin{definition}
\label{def2.7} Let $X$ be a smooth projective variety over $k$ and
let $c_1,\dots ,c_n$ be $0$-cycles on $X$. We define their {\it
wedge product} to be the following:
   $$
   c_1\wedge \dots \wedge c_n=\frac{1}{n!}
   \sum _{\sigma \in \Sigma _n}
   sgn(\sigma )c_{\sigma (1)}\times \dots \times c_{\sigma (n)}
   $$
where $c_{\sigma (1)}\times \dots \times c_{\sigma (n)}$ is the
exterior product of cycles, see \cite{Ful}, Ch.5.
\end{definition}

In \cite{Kim}, 5.14, it is proved that, if a surface $X$ is the
product of $2$ curves, then there exists an integer $N$, such that
the product $c_1\wedge \dots \wedge c_N=0$, where $c_i$ are
$0$-cycles in the Albanese kernel $T(X)$. The following theorem
extends this result to any surface $X$ whose motive is finite
dimensional.

\begin{theorem}
\label{Prop8} Let $X$ be a smooth projective surface over $k$. If
the motive $M(X)\in \mathcal M_{rat}$ is finite dimensional then
$c_1\wedge \dots \wedge c_{d+1}=0$, where $c_i$ are $0$-cycles in
the Albanese kernel $T(X)$, $d=b_2-\rho $, $b_2=dim(H^2(X))$ and
$\rho =dim(NS(X)_{\mathbb Q})$.
\end{theorem}

\begin{pf} The motive $M(X)$ has a Chow-K\"unneth
decomposition as follows:
   $$
   M(X)=\sum _{0\le i\le 4}M_i(X)\; .
   $$
Since $M(X)$ is finite dimensional, $M_2(X)$ is also finite
dimensional. From \cite{Mu3} it follows that $A^1(M_2(X))
=NS(X)_{\mathbb Q}$ and from \cite{Sch}, 2.2:
   $$
   Hom_{\mathcal M_{rat}}(M_2(X),\mathbb L)\simeq
   Hom_{\mathcal M_{rat}}(\mathbb L,M_2(X))\simeq
   NS(X)_{\mathbb Q}\; ,
   $$
where $\mathbb L =(Spec(k),id,-1)$ is the Lefschetz motive.
$NS(X)_{\mathbb Q}$ is a finite dimensional $\mathbb Q$-vector
space of dimension $\rho $ (the corank of $Pic(X)_{\mathbb Q}$).
Let $[e_i]$, for $1\le i\le \rho $, be a base for $NS(X)_{\mathbb
Q}$ and let $\alpha =\sum q_i[e_i]\in NS(X)_{\mathbb Q}$. Let
$f_{\alpha }:\mathbb L\to M_2(X)$ be the corresponding morphism in
$\mathcal M_{rat}$. Then $f_{\alpha }=\sum q_i[Spec(k)\times
e_i]$. The transpose $f^t_{\alpha }$ is a morphism $M_2(X)\to
\mathbb L$ and
   $$
   f^t_{\alpha }\circ f_{\alpha }\in
   Hom_{\mathcal M_{rat}}(\mathbb L,\mathbb L)
   \simeq \mathbb Q\; .
   $$
Therefore, for every $i\le \rho $, $f_{[e_i]}:\mathbb L\to M_2(X)$
is an injective map. Let $f= \sum f_{[e_i]}$. Then $f$ defines an
injective map:
   $$
   \mathbb L\oplus \dots \oplus
   \mathbb L \; (\rho \; \hbox{times})\to M_2(X)\; .
   $$

This yields a splitting in $\mathcal M_{rat}$:
   $$
   M_2(X)=\rho \mathbb L\oplus N\; .
   $$

We have: $H^i(N)=0$ for $i\ne 2$ and $H^2(M_2(X))=H^2(X)= \rho
H^2(\mathbb L)\oplus H^2(N)$ where $H^2(\mathbb L)=\mathbb Q$.
Therefore, $H^2(N)=(b_2-\rho )\cdot \mathbb Q$.

$M_2(X)$ being finite dimensional $N$ is finite dimensional too.
$N$ is evenly finite dimensional because it does not have any odd
cohomology, see \cite{Kim}, 3.9. Therefore,
$dim(N)=dim(H^2(N))=kim(N)= d=(b_2-\rho )$, so that $\wedge
^{d+1}N=0$. We also have $A^2(N)=A^2(M_2(X))=T(X)$. If
$N=(Y,q,n)\in \mathcal M_{rat}$, then it follows from the
definition of $\wedge $ that, for any $r$, $\wedge ^rN$ is the
image of the motive $N^r$ under the projector $(1/r!)(\sum
_{\sigma \in \Sigma _r}sgn(\sigma ) \Gamma _{\sigma }\circ q^r)$.
Therefore, if $c_1,\dots ,c_{d+1}$ are $0$-cycles in $T(X)$, then
the cycle $c=c_1\wedge \dots \wedge c_{d+1}$ belongs to
$A^2(\wedge ^{d+1}N)=0$. This proves that $c=0$.
\end{pf}

\section{Surfaces with $p_g =0$}
\label{S3}

 From Th. \ref{Prop8} it follows that, if $X$ is a smooth
projective surface with $p_g=0$ (a condition which is equivalent
to $b_2=\rho $), then the finite dimensionality of the motive
$M(X)$ implies $T(X)=0$. In this section we prove (Theorem
\ref{th12}) that the converse also holds.

We first recall, see \cite{BV}, the definition of a {\it balanced}
variety:

\begin{definition}
\label{def3.1} Let $X$ be a reduced separated and equidimensional
scheme of finite type over $k$ and $d=dim(X)$. $X$ is said to be
balanced of weight $w$ if there exist cycles $\Gamma _1$ and
$\Gamma _2$ of codimension $d$ on $X\times X$, such that
   $$
   [\Delta _X]=[\Gamma _1]+[\Gamma _2]
   $$
in $A^d(X\times X)$, where $\Gamma_1$ is supported on $Z_1\times
X$, $\Gamma _2$ is supported on $X\times Z_2$, $Z_1$ and $Z_2$ are
equidimensional closed subschemes of $X$ and
   $$
   w=min\{ dim(Z_1),dim(Z_2)\} \; .
   $$
\end{definition}

\begin{lemma}
\label{lemma9} Let $\mathcal C$ be a pseudoabelian category and
let $X$ be an object in $\mathcal C$. Assume that there exists a
finite collection of objects $Y_i$ and morphisms
$X\stackrel{a_i}{\to }Y_i\stackrel{b_i}{\to }X$ in $\mathcal C$,
$i=1,\dots ,n$, such that $id_X=\sum _{i=1}^nb_ia_i$ in
$End_{\mathcal C}(X)$. Let
     $$
     X\stackrel{f}{\longrightarrow }\bigoplus _{i=1}^nY_i
     \stackrel{g}{\longrightarrow }X
     $$
be the morphisms induced by $\{ a_i\} $ and $\{ b_i\} $
respectively. Then $gf=1_X$ and therefore $X$ is isomorphic to a
direct summand of $\oplus _{i=1}^nY_i$.
\end{lemma}

\begin{pf}
Let $\pi _i:\oplus _{i=1}^nY_i\to Y_i$ and $\iota _i:Y_i\to \oplus
_{i=1}^nY_i$ be the canonical projection and the canonical
embedding. Then
   $$
   gf=g\circ 1_{\oplus _{i=1}^nY_i}\circ f=
   g\circ \left(\sum _{i=1}^n\iota _i\pi _i\right)\circ f=
   \sum _{i=1}^ng\iota _i\pi _if=\sum _{i=1}^nb_ia_i=1_X\; .
   $$
\end{pf}

\begin{theorem}
\label{th10} Let $X$ be a smooth projective (equidimensional)
variety balanced by subschemes $Z_1$ and $Z_2$. Let $\tilde Z_i$
be a desingularization of $Z_i$, $i=1,2$. Assume that the motives
$M(\tilde Z_1)$ and $M(\tilde Z_2)$ are finite dimensional. Then
the motive $M(X)$ is finite dimensional.
\end{theorem}

\begin{pf} For any $Y\in \mathcal V_k$ and $\alpha \in
Corr^0(Y,X)$ let $f_{\alpha }:M(Y)\to M(X)$ be the corresponding
morphism from $M(Y)$ to $M(X)$. By assumptions, $[\Delta
_X]=[\Gamma _1]+[\Gamma _2]$ in $Corr^0(X,X)$, where $\Gamma
_1\subset Z_1\times X$ and $\Gamma _2\subset X\times Z_2$. In
other words, $id_{M(X)}=f_{[\Gamma _1]}+f_{[\Gamma _2]}$ in
$End_{\mathcal M_{rat}}(M(X))$.

Let $s_i:\tilde Z_1\to X$ be a composition of the closed embedding
$Z_i\hookrightarrow X$ with a blow up desingularization
$v_i:\tilde Z_i\to Z_i$ of $Z_i$, $i=1,2$. Let also
$w_i=dim(Z_i)$.

Since $\Gamma _1$ lies on $Z_1\times X$, we may consider its class
$[\Gamma _1]$ in the Chow group $A^{w_1}(Z_1\times X)$ of the
scheme $Z_1\times X$. Let $[\tilde \Gamma _1]$ be the pull back of
$[\Gamma _1]$ with respect to the morphism $v_1\times id_X:\tilde
Z_1\times X\to Z_1\times X$, see \cite{Ful}. The variety $\tilde
Z_1$ is smooth projective and $dim(\tilde Z_1)=w_1$. It follows
that $[\tilde \Gamma _1]$ lies in $A^{w_1}(\tilde Z_1\times
X)=Corr^0(\tilde Z_1,X)$. Consider the corresponding morphism
$f_{[\tilde \Gamma _1]}:M(\tilde Z_1)\to M(X)$ in the category
$\mathcal M_{rat}$. Since $v_1\times id_X$ is a blow up, it
follows that $(v_1\times id_X)_*(v_1\times id_X)^*([\Gamma
_1])=[\Gamma _1]$ \cite{Ful}. Therefore we get:
$f_{[\Gamma_1]}=f_{[\tilde \Gamma_1]}\circ M(s_1)$, whence the
morphism $f_{[\Gamma_1]}:M(X)\to M(X)$ factors through the motive
$M(\tilde Z_1)$.

Similarly one shows, by applying duality in $\mathcal M_{rat}$,
that the morphism $f_{[\Gamma _2]}:M(X)\to M(X)$ factors through
the motive $M(\tilde Z_2)(w_2-d)$ where $d=dim (X)$. Indeed, let
$\Gamma _2^t$ be the transposition of the cycle $\Gamma _2$. Let
$[\tilde \Gamma _2^t]$ be a pull back of its class $[\Gamma _2^t]$
(in the Chow group $A^{w_2}(Z_2\times X)$) with respect to the
blow up $v_2\times id_X:\tilde Z_2\times X\to Z_2\times X$. As
above we get: $f_{[\Gamma _2^t]}=f_{[\tilde \Gamma _2^t]}\circ
M(j_2)$. Considering $f_{[\Gamma _2^t]}$ as an endomorphism of the
motive $M(X)(d)=M(X)\otimes \mathbb L^{-d}$ it factors through
$M(\tilde Z_2)(d)$. By dualizing we see that $f_{[\Gamma
_2]}:M(X)\to M(X)$ factors through the motive $M(\tilde
Z_2)(w_2-d)$.

By Lemma \ref{lemma9} we have that $M(X)$ is isomorphic to a
direct summand of the motive $M(\tilde Z_1)\oplus M(\tilde
Z_2)(w_2-d)$. Since both motives $M(\tilde Z_1)$ and $M(\tilde
Z_2)(w_2-d)$ are finite dimensional, their direct sum is finite
dimensional. Therefore $M(X)$ is also finite dimensional.
\end{pf}

\begin{corollary}
\label{corollary11} Let $X$ a smooth projective surface. Assume
that either $X$ or a Zariski open dense subset $U$ of $X$ are
balanced. Then $M(X)$ is finite dimensional.
\end{corollary}

\begin{pf}
Let $Z=X-U$. Then $codim_X(Z)\le 1$. Let $U$ be balanced over
closed subschemes $Z_1$ and $Z_2$ of codimension $\le 1$. By a
result of Barbieri-Viale, \cite{BV}, $X$ is balanced of weight
$\le 1$. The motives of points and curves are finite dimensional.
Therefore Theorem \ref{th10} implies that $M(X)$ is finite
dimensional.
\end{pf}
\medskip

\begin{remark}
{\rm For any field $k$ of characteristic $0$ V.Voevodsky has
constructed in \cite{Voev} a triangulated category of motives
$DM(k)$ over $k$ and a functor $M:Sm/k\to DM(k)$ from the category
$Sm/k$ of smooth separated schemes over $k$ into $DM(k)$. This
triangulated category contains a full subcategory, generated by
motives $M(X)$ of smooth projective varieties $X$, which is
equivalent to $\mathcal M_{rat}$. Moreover, it is pseudo-abelian
and, if we consider finite correspondences on schemes with
coefficients in $\mathbb Q$ to construct $DM(k)$, it is $\mathbb
Q$-linear. Therefore we can define, according to Def.
\ref{Def2.1}, finite dimensionality of the motive $M(V)$ for every
$V\in Sm/k$. Moreover, if $U$ is an open subset of a smooth
projective variety $X$ one has the following distinguished
triangle in $DM(k)$ \cite{Voev}, 3.5.4:
   $$
   M(U)\to M(X)\to M(Z)(i)[2i]\to M(U)[1]
   $$
where $Z=X-U$ and $i$ is the codimension of $Z$ in $X$. If $X$ is
a surface and $U$ an open subset of $X$ then $Z$ has dimension
$\le 1$, so that $M(Z)$ is finite dimensional. This implies that
also $M(Z)(i)[2i]$ is finite dimensional. Therefore Corollary
\ref{corollary11} naturally suggests the following question:
assuming that $M(U)$ is finite dimensional, is $M(X)$ also finite
dimensional?}
\end{remark}

The next result (Theorem \ref{th12}) shows that, for a surface
with $p_g=0$ finite dimensionality of the motive $M(X)$ is
equivalent to the finite dimensionality of the Chow group of
$0$-cycles in the sense of Mumford. Here is the definition:

\begin{definition}
\label{def3.3} Let $X$ a smooth projective variety of dimension
$d$ over an algebraically closed field $k$ and let $A^d(X)_0$ be
the group of $0$-cycles of degree $0$ on $X$. Then $A^d(X)_0$ is
{\it finite dimensional} if there exists an integer $n$, such that
the natural map
   $$
   s_n:S^nX\times S^nX\to A^d(X)_0
   $$
is surjective, where $s_n(A,B)=A-B$ and $S^nX$ is the $n$-th
symmetric power of $X$.
\end{definition}

\begin{remark}
{\rm Note that, if $p_g>0$, then finite dimensionality of the
motive $M(X)$ does not, in general, imply the finite
dimensionality of the Chow group, as it can be shown by taking
products of curves $C_i$ of genus $>1$. If $X$ is a complex
surface with $p_g=0$, then $A^2(X)_0$ is finite dimensional iff
Conjecture \ref{bloch2} holds for $X$.}
\end{remark}

\begin{theorem}
\label{th12} Let $X$ be a smooth projective surface over an
algebraically closed field $k$ of characteristic $0$ with
$p_g(X)=0$. Then the motive $M(X)$ is finite dimensional if and
only if the group $A^2(X)_0$ is finite dimensional (i.e. Bloch's
conjecture on Albanese kernel is true for $X$).
\end{theorem}

\begin{pf}
If $M(X)$ is finite dimensional then by \cite{GP1}, Theorem 2.11,
we have $ker(\pi _3)=T(X)=0$ where $T(X)$ is the Albanese kernel.
This implies that $A^2(X)_0$ is finite dimensional, see \cite{J1},
1.6.

Conversely, assume that $A^2(X)_0$ is finite dimensional. Then
there exists, \cite{J1}, 1.6, a closed subscheme $Y$ of dimension
$\le 1$, such that $A^2(X-Y)=0$. By results of \cite{BS} $X$ is
balanced of weight $\le 1$. From Theorem \ref{th10} it follows
that $M(X) $ is finite dimensional.
\end{pf}

\begin{remark}
(Relations with K-theory) {\rm Let $X$ be a smooth projective
surface over $\mathbb C$. Then one has the following description
for the $K$-groups $K_n(X)$, for $n>0$, see \cite{PW1}, 6.7:
    $$
    K_n(X)\simeq \left\{
    \begin{array}{ll}
    B\oplus (\mathbb Q/\mathbb Z)^{2+b_2}\oplus V_n & \mbox{if $n\geq 1$ odd} \\
    A\oplus (\mathbb Q/\mathbb Z)^{b_1+b_3}\oplus V_n & \mbox{if $n\geq 2$ even}
    \end{array}
    \right.
    $$
where $A=H^2(X(\mathbb C),\mathbb Z)_{tors}$, $B=H^3(X(\mathbb
C),\mathbb Z)_{tors}$, $b_i$ are the Betti numbers and $V_n$ are
uniquely divisible groups.

A similar result also holds for any smooth variety over $\mathbb
C$ \cite{PW2} if one assumes the so called {\it norm residue
Conjecture} which asserts that the norm residue map: $K^M_n(F)/m
\to H^n_{et}(F,\mathbb Z/m)$ is an isomorphism for all $m$, where
$F$ is the function field of $X$ and $K^M_*$ is Milnor's
$K$-theory.

It follows that, for a surface $X$, $K_n(X)_{tors}$ depends only
upon the topological invariants of the manifold $X(\mathbb C)$. On
the other hand the groups $K_n(X)_{\mathbb Q}=K_n(X)\otimes
\mathbb Q$ depend on the motive $M(X)$ via the Bloch-Lichtenbaum
spectral sequence $E^{p,q}_2=H^{p-q}_{\mathcal M}(X,\mathbb
Q(-q))$ which converges to $K_{-p-q}(X)_{\mathbb Q}$. Here
   $$
   H^{2i}_{\mathcal M}(X,\mathbb Q(i))=
   Hom_{DM(k)}(M(X),\mathbb Q(i)[2i])
   $$
is the motivic cohomology of $X$ and $\mathbb Q(i)[2i]$ plays a
role of the power $\mathbb L^i$ in $DM(k)$ \cite{Voev}.

Now let $X$ be a smooth projective surface with $p_g=q=0$. If
$M(X)$ is finite dimensional then, by \cite{GP1}, 2.14, the motive
$M(X)$ is "trivial" in the sense that it is a direct sum of the
unit motive 1, of $\mathbb L^2$ and of a finite number of copies
of $\mathbb L$. From Th. \ref{th12} it follows that the Albanese
kernel $T(X)$ vanishes and this, by \cite{Pe}, Th.0.1, implies
   $$
   K_n(X)_{\mathbb Q}\simeq
   K_0(X)\otimes K_n(\mathbb C)_{\mathbb Q}
   $$
So also the higher $K$-theory of $X$ is "trivial".

Note that, if either $p_g$ or $q$ do not vanish, then the above
isomorphism is, in general, not true, see \cite{Pe}.}
\end{remark}

% The endnotes section will be placed here.

%\theendnotes

\begin{small}

\end{small}

\bigskip

\noindent {\tt guletskii@im.bas-net.by}

\noindent {\sc Institute of Mathematics, Surganova 11, 220072
Minsk, Belarus}

\medskip

\noindent {\tt pedrini@dima.unige.it}

\noindent {\sc Dipartimento di Matematica, Universit\`a di Genova,
Via Dodecaneso 35, 16146 Genova, Italy}

\end{document}